\documentclass[a4paper,11pt]{amsart}

\usepackage[utf8]{inputenc}
\usepackage[T1]{fontenc}
\usepackage{lmodern}
\usepackage{microtype}
\usepackage{amsmath,amssymb,mathtools}
\usepackage{enumitem}
\usepackage{booktabs,tabularx}
\usepackage{placeins}
\usepackage[hidelinks]{hyperref}
\hypersetup{pdftitle={Rigid Lie Algebras from Semisimple Lie Algebras and Their Symplectic Modules},
pdfauthor={Estela Fernandez and Paulo Tirao}}
\newtheorem{theorem}{Theorem}[section]
\newtheorem{proposition}[theorem]{Proposition}
\newtheorem{lemma}[theorem]{Lemma}
\newtheorem{corollary}[theorem]{Corollary}

\theoremstyle{definition}
\newtheorem{definition}[theorem]{Definition}
\newtheorem{example}[theorem]{Example}

\theoremstyle{remark}
\newtheorem{remark}[theorem]{Remark}

\newcommand{\C}{\mathbb{C}}
\newcommand{\g}{\mathfrak{g}}
\newcommand{\h}{\mathfrak{h}}
\newcommand{\s}{\mathfrak{s}}
\newcommand{\Der}{\operatorname{Der}}
\newcommand{\Hom}{\operatorname{Hom}}

\title[Rigid Lie Algebras from Symplectic Modules]
{Rigid Lie Algebras from Semisimple Lie Algebras and Their
Symplectic Modules}

\author{Estela Fern\'andez}
\address{FaCET, Universidad Nacional de Tucum\'an, Argentina}
\email{efernandez@herrera.unt.edu.ar}

\author{Paulo Tirao}
\address{CIEM -- FaMAF, CONICET -- Universidad Nacional de C\'ordoba, Argentina}
\address{Guangdong Technion Israel Institute of Technology, China}
\email{paulotirao@unc.edu.ar}


\date{}

\begin{document}

\begin{abstract}
We consider two extensions associated with a nonzero finite-dimensional
symplectic module $V$ over a complex semisimple Lie algebra
$\mathfrak{s}$, and study their rigidity. The toral extension adjoins
a derivation acting as the identity on $V$ and trivially on
$\mathfrak{s}$; the Heisenberg extension uses the invariant symplectic
form on $V$ to define a central bracket.

For the toral extension, we compute the first and second adjoint
cohomology groups in terms of $\operatorname{End}_{\mathfrak{s}}(V)$,
and obtain rigidity precisely when $V$ is irreducible. These results
hold for arbitrary nonzero finite-dimensional modules.
For the Heisenberg extension, we express the second adjoint cohomology
as the kernel of a map combining the symplectic centralizer with an
equivariant Jacobi map. A Killing-form argument shows that every
multiplicity-free sum of irreducible symplectic modules gives an
algebraically rigid Heisenberg extension. We also treat constructions
with independent toral generators or Heisenberg centers. Repeated
summands lead to explicit nonrigid examples, including a deformation
to $\mathfrak{sl}_3$. The proofs use invariant cochains and recover
several previously studied families.\end{abstract}

\maketitle

\section{Introduction}

All Lie algebras and modules in our results are finite-dimensional over
$\C$. A Lie algebra is \emph{geometrically rigid} if its change-of-basis
orbit is Zariski open in the variety of Lie algebra laws, and
\emph{algebraically rigid} if its second adjoint cohomology vanishes.
Algebraic rigidity implies geometric rigidity
\cite{NijenhuisRichardson1967}; the converse fails in general
\cite{Richardson1967}.

Let $\s$ be semisimple and let $(V,\omega)$ be a nonzero symplectic
$\s$-module. Regard $V$ as abelian. The direct product
\begin{equation}\label{eq:intro-starting-algebra}
\g_0=(\s\ltimes V)\times\C a
\end{equation}
admits two linear deformations, defined by
\[
\tau(a,v)=v,\qquad \eta(u,v)=\omega(u,v)a,
\]
with all unspecified values zero. Their nonzero members have respective
isomorphism types
\[
\g_\tau=(\s\oplus\C a)\ltimes V,
\qquad
\g_\eta=\s\ltimes\h(V,\omega).
\]
In the first algebra, $a$ acts as the identity on $V$. In the second,
$\h(V,\omega)=V\oplus\C a$ has bracket $[u,v]=\omega(u,v)a$.
We compute the low-degree adjoint cohomology of these targets and show
how their rigidity depends on the irreducible decomposition of $V$.

The toral calculation applies to any nonzero $\s$-module, without a
symplectic assumption. With $A=\operatorname{End}_\s(V)$, it gives
\[
H^0(\g_\tau,\g_\tau)=0,\qquad
H^1(\g_\tau,\g_\tau)\simeq H^2(\g_\tau,\g_\tau)
\simeq A/\C\operatorname{Id}_V
\]
(Proposition~\ref{prop:reducible-tau}). Thus $\g_\tau$ is complete and
algebraically rigid when $V$ is irreducible. When $V$ is reducible,
equivariant projections give nontrivial linear deformations, proving
geometric nonrigidity. Adjoining an independent toral generator for each
irreducible component instead produces an extension whose adjoint
cohomology vanishes in every degree (Proposition~\ref{prop:full-torus}).

Our main Heisenberg result applies to every nonzero multiplicity-free
sum of irreducible symplectic $\s$-modules:
\begin{equation}\label{eq:intro-eta-rigidity}
H^0(\g_\eta,\g_\eta)\simeq\C,\qquad
H^1(\g_\eta,\g_\eta)\simeq\C,\qquad
H^2(\g_\eta,\g_\eta)=0.
\end{equation}
No exclusion of adjoint components from $\Lambda^2V$ is required
(Theorem~\ref{thm:multiplicity-free-eta}). The proof combines an invariant
cochain calculation with a Killing-form argument: an equivariant map
$\Lambda^2V\to\s$ satisfying the corresponding Jacobi identity must
vanish. The same argument proves rigidity when each irreducible
symplectic component has its own Heisenberg center
(Theorem~\ref{thm:multi-heisenberg}). These algebras have a single
diagonally acting Levi factor.

The general calculation underlying \eqref{eq:intro-eta-rigidity}
expresses $H^2(\g_\eta,\g_\eta)$ as the kernel of an explicit map
involving the symplectic centralizer of the $\s$-action and the
adjoint-valued two-cochains (Proposition~\ref{prop:reducible-eta}). This
formula retains the cancellations that occur when irreducible
components are repeated. We exhibit two resulting failures of rigidity:
three or more isomorphic Heisenberg blocks admit explicit central-valued
deformations, while two copies of the defining $\mathfrak{sl}_2$-module
with one shared center admit a deformation to $\mathfrak{sl}_3$.

These results belong to the cohomological approach to rigidity developed
from the Hochschild--Serre factorization theorem and equivariant
cochains. Richardson's examples distinguish geometric from algebraic
rigidity for semidirect products with abelian radical
\cite{Richardson1967}. Related semidirect products involving
$\mathfrak{so}(3)$ were studied in \cite{Campoamor2021}.
Burde and Wagemann give explicit adjoint-cohomology calculations for
semidirect products and low-dimensional perfect Lie algebras, including
examples with Heisenberg radical
\cite{BurdeWagemann2023,BurdeWagemann2024}.

The toral extensions fit the external-torus framework of Carles
\cite{Carles1984,Carles1985}; his cohomological study of generalizations
of Borel subalgebras is also relevant to the full-torus calculation
\cite{Carles1993}. We give direct proofs for the extensions considered
here and separate them from these broader structural results. The
deformation defined by $\tau$ is the split codimension-one construction
of Grunewald and O'Halloran \cite{GrunewaldOHalloran1993}; see also
\cite{HerreraGranadaTirao2016}. The deformation defined by $\eta$ instead
changes the bracket on $\Lambda^2V$.

Several irreducible Heisenberg families are already known.
Campoamor-Stursberg proved cohomological rigidity for the Schr\"odinger
algebra and its mass central extension in the nonexceptional spatial
dimensions \cite{Campoamor2017}. Khudoyberdiyev and Jumaniyozov obtained
the corresponding vanishing results for mass central extensions of
conformal Galilei algebras \cite[Theorems~5.1 and~5.2]{KhudoyberdiyevJumaniyozov2025}.
Their $d=1$ case includes all Heisenberg extensions by odd irreducible
$\mathfrak{sl}_2$-modules. We recover these results as applications of
the general theorem. The symplectic oscillator algebras
$\mathfrak{sp}_{2n}\ltimes\h_{2n+1}$ provide another family
\cite{LiuZhao2019}. Complete extensions with Heisenberg nilradical,
which include additional derivations, were studied in
\cite{Meng1994,JiangMeng1998}.

The emphasis of this paper is the general multiplicity-free rigidity
criterion, the kernel formula that extends beyond it, and the comparison
between a shared center and separate Heisenberg centers. The explicit
deformations describe some of the limitations of these criteria.
Section~\ref{sec:preliminaries} collects the cohomological and
representation-theoretic tools and constructs the common source.
Section~\ref{sec:toral} treats toral extensions.
Section~\ref{sec:heisenberg} proves the Heisenberg results, and
Section~\ref{sec:examples} gives applications and nonrigid examples.

\subsection*{Acknowledgments}
Some of the results of this paper were obtained as part of the Ph.D. thesis
of the first author. The authors thank Edison Fern\'andez-Culma and Nadina
Rojas for valuable discussions concerning that thesis.

\section{Cohomology, symplectic modules, and the common source}
\label{sec:preliminaries}

We use $\g$ for a general Lie algebra, $\s$ for a complex semisimple Lie
algebra, and $V$ for an $\s$-module. The symbol $\g_0$ denotes the
undeformed direct product in \eqref{eq:intro-starting-algebra}; its
targets are $\g_\tau$ and $\g_\eta$. A superscript $(r)$ indicates a
construction with $r$ added generators. The corresponding radicals are
denoted by $\mathfrak r_0$, $\mathfrak r_\tau$, and $\mathfrak r_\eta$,
with the same superscript convention. We write $\C a$ for the line
spanned by $a$, and $\C$ for the scalar field or the trivial module.
\subsection{Geometric and algebraic rigidity}

Let $\g$ be a finite-dimensional complex Lie algebra with bracket $\mu$,
and let $\mathcal{L}(\g)$ denote the variety of Lie algebra laws on its
underlying vector space. The group $\operatorname{GL}(\g)$ acts on
$\mathcal{L}(\g)$ by change of basis.

\begin{definition}
The Lie algebra $\g$ is \emph{geometrically rigid} if the orbit of
$\mu$ under $\operatorname{GL}(\g)$ is Zariski open in
$\mathcal{L}(\g)$. It is \emph{algebraically rigid}
(also called \emph{cohomologically rigid}) if
\[
H^2(\g,\g)=0,
\]
where the coefficients are taken in the adjoint module.
\end{definition}

Algebraic rigidity implies geometric rigidity
\cite{NijenhuisRichardson1967}. The converse does not hold in general
\cite{Richardson1967}.

\begin{definition}
A Lie algebra $\g$ is \emph{complete} if its center is zero and every
derivation is inner; equivalently,
\[
H^0(\g,\g)=H^1(\g,\g)=0.
\]
\end{definition}

\subsection{Cohomology of semidirect products}

For a Lie algebra $\g$ and a $\g$-module $M$, write
$C^q(\g,M)=\Hom(\Lambda^q\g,M)$. Our differential is
\begin{align*}
(df)(x_0,\ldots,x_q)
={}&\sum_{i=0}^q(-1)^i x_i\cdot
 f(x_0,\ldots,\widehat{x_i},\ldots,x_q)\\
&+\sum_{0\leq i<j\leq q}(-1)^{i+j}
 f([x_i,x_j],x_0,\ldots,\widehat{x_i},\ldots,
 \widehat{x_j},\ldots,x_q).
\end{align*}
In particular, for adjoint coefficients, $(dz)(x)=[x,z]$ and
$(df)(x,y)=[x,f(y)]-[y,f(x)]-f([x,y])$.
All displayed cohomology isomorphisms are vector-space isomorphisms
unless explicitly stated otherwise.

Let $I$ be an ideal of $\g$. The actions of $\g$ on $I$ and $M$
induce an action on $C^\bullet(I,M)$, given by
\[
(x\cdot f)(u_1,\ldots,u_q)
=x\cdot f(u_1,\ldots,u_q)
-\sum_{j=1}^q f(u_1,\ldots,[x,u_j],\ldots,u_q).
\]
This action commutes with the differential and therefore induces a
$\g$-module structure on $H^q(I,M)$.

We use the following form of the Hochschild--Serre factorization theorem
\cite[Theorem~13, p.~603]{HochschildSerre1953}.

\begin{theorem}[Hochschild--Serre factorization theorem]
\label{thm:HS-factorization}
Let $\g$ be a finite-dimensional complex Lie algebra, let $M$ be a
finite-dimensional $\g$-module, and let $I$ be an ideal such that
$\g/I$ is semisimple. Then, for every $q\geq0$,
\begin{equation}\label{eq:HS-factorization}
H^q(\g,M)
\simeq
\sum_{i+j=q}H^i(\g/I,\C)\otimes H^j(I,M)^\g.
\end{equation}
Here the summation denotes a direct sum of vector spaces.
\end{theorem}

Since $\g/I$ is semisimple, we have $H^0(\g/I,\C)=\C$ and,
by Whitehead's lemmas,
\[
H^1(\g/I,\C)=H^2(\g/I,\C)=0.
\]
Consequently, restriction induces isomorphisms
\begin{equation}\label{eq:HS-low-degree-g}
H^q(\g,M)\simeq H^q(I,M)^\g,
\qquad q=0,1,2.
\end{equation}

Now write $\g=\s\ltimes I$, where $\s$ is a semisimple subalgebra
complementary to $I$. By \cite[(1.2), p.~1011]{Leger1957},
\begin{equation}\label{eq:cohomology-invariants}
H^q(I,M)^\g=H^q(I,M)^\s
\qquad(q\geq0).
\end{equation}
Indeed, the action of $I$ on its cohomology is trivial. Combining
\eqref{eq:HS-low-degree-g} and \eqref{eq:cohomology-invariants} gives
the reduction used below.

\begin{corollary}\label{cor:HS-low-degree}
Let $\g=\s\ltimes I$, with $\s$ semisimple, and let $M$ be a
finite-dimensional $\g$-module. Then, for $q=0,1,2$,
\begin{equation}\label{eq:HS-low-degree}
H^q(\g,M)\simeq H^q(I,M)^\s.
\end{equation}
\end{corollary}

Moreover, every $\s$-invariant cohomology class has an $\s$-invariant
cocycle representative, and every $\s$-invariant coboundary has an
$\s$-invariant primitive \cite[(1.3)--(1.4), p.~1011]{Leger1957}.
Equivalently, taking $\s$-invariants is exact on finite-dimensional
$\s$-modules, and
\[
\bigl(B^q(I,M)\bigr)^\s=d\Bigl(\bigl(C^{q-1}(I,M)\bigr)^\s\Bigr).
\]
Thus
\begin{equation}\label{eq:invariant-cohomology}
H^q(I,M)^\s
\simeq
\frac{\bigl(Z^q(I,M)\bigr)^\s}{d\Bigl(\bigl(C^{q-1}(I,M)\bigr)^\s\Bigr)}
=
H^q\bigl(\bigl(C^\bullet(I,M)\bigr)^\s\bigr),
\end{equation}
where $C^{-1}(I,M)=0$. The invariant complex in the last expression
has degree-$q$ term $\Hom_\s(\Lambda^q I,M)$, with the restriction
of the Chevalley--Eilenberg differential. Invariance is imposed on
cochains, including both their arguments and their values.

These reductions apply to both abelian and nonabelian radicals. In each
case the differential uses the bracket and coefficient action of the
algebra being considered.

\subsection{Zero-weight cochains}
We also use the following consequence of the Cartan homotopy formula.  Let $\g$ be a Lie algebra, let $M$ be
a $\g$-module, and let $x\in\g$.  Denote by
$\mathcal L_x$ the endomorphism induced by $x$ on the cochain complex
$C^\bullet(\g,M)$ and put
$(\iota_xf)(y_1,\ldots,y_{q-1})=f(x,y_1,\ldots,y_{q-1})$.
Then
\begin{equation}\label{eq:Cartan-homotopy}
\mathcal L_x=d\iota_x+\iota_xd.
\end{equation}
Here
\[
(\mathcal L_xf)(y_1,\ldots,y_q)
=x\cdot f(y_1,\ldots,y_q)
-\sum_{i=1}^qf(y_1,\ldots,[x,y_i],\ldots,y_q).
\]
Since $d\mathcal L_x=\mathcal L_xd$, the differential preserves every
eigenspace of $\mathcal L_x$. When $\mathcal L_x$ is diagonalizable,
these eigenspaces therefore give a direct-sum decomposition of the
cochain complex into subcomplexes. If $f$ is a cocycle of eigenvalue
$\lambda\neq0$, Cartan's identity gives
\[
f=d\bigl(\lambda^{-1}\iota_x f\bigr).
\]
Thus every nonzero-eigenvalue subcomplex has zero cohomology, and the
inclusion of the zero-eigenvalue (zero-weight) subcomplex induces an
isomorphism in cohomology. The same argument applies to
$\bigl(C^\bullet(I,M)\bigr)^\s$ when $x\in I$ and $[x,\s]=0$, since insertion of
$x$ then preserves $\s$-invariance.
\subsection{Symplectic modules}\label{sec:symplectic-modules}

Let $\s$ be a semisimple Lie algebra and let $V$ be an $\s$-module.

\begin{definition}
The module $V$ is called \emph{symplectic} if it admits a nondegenerate
$\s$-invariant alternating bilinear form
\[
\omega\colon V\times V\longrightarrow\C.
\]
Thus
\[
\omega(x\cdot v,v')+\omega(v,x\cdot v')=0
\]
for every $x\in\s$ and $v,v'\in V$.
\end{definition}

If $V$ is irreducible, such a form, when it exists, is unique up to a nonzero
scalar.  Indeed, an invariant bilinear form is the same as an
$\s$-homomorphism $V\to V^*$.  Hence an irreducible module admits a
nonzero invariant bilinear form if and only if it is self-dual, and any such
form is automatically nondegenerate.  By Schur's lemma, it is then either
symmetric or alternating.

This distinction is measured by the Frobenius--Schur indicator.  If
$c$ denotes the transposition
\[
c(v\otimes v')=v'\otimes v,
\]
we denote the indicator of an irreducible $\s$-module $V$ by
\begin{equation}\label{eq:FS-indicator}
\nu(V)=\operatorname{tr}
\left(c\bigm|(V\otimes V)^\s\right).
\end{equation}
The Frobenius--Schur theorem gives
\begin{equation}\label{eq:FS-types}
\nu(V)=
\begin{cases}
 0,  & V\not\simeq V^*,\\
 1,  & V\text{ admits an invariant symmetric form},\\
-1,  & V\text{ admits an invariant alternating form}.
\end{cases}
\end{equation}
See \cite{AbuHamedGelaki2007} for the Frobenius--Schur indicators of
finite-dimensional modules over complex semisimple Lie algebras.  In
particular, for an irreducible module,
\begin{equation}\label{eq:symplectic-FS}
V\text{ is symplectic}
\quad\Longleftrightarrow\quad
\nu(V)=-1
\quad\Longleftrightarrow\quad
\Hom_\s(\Lambda^2V,\C)\neq0.
\end{equation}
In this case the last space is one-dimensional and $\dim V$ is even.

The indicator is also computable directly from the highest weight.  Let
$V(\lambda)$ be the irreducible module with highest weight $\lambda$, let
$w_0$ be the longest element of the Weyl group, and let $\rho^\vee$ be the
half-sum of the positive coroots.  Tits's formula \cite[Theorem~3.6]{AbuHamedGelaki2007} states that
\begin{equation}\label{eq:Tits-indicator}
\nu(V(\lambda))=
\begin{cases}
0, & \lambda\neq-w_0\lambda,\\
(-1)^{\lambda(2\rho^\vee)}, & \lambda=-w_0\lambda.
\end{cases}
\end{equation}
Thus the symplectic condition is an explicit parity condition on a self-dual
highest weight.

For a general module, write
\[
V=\bigoplus_{i=1}^k U_i\otimes M_i,
\]
where the $U_i$ are pairwise nonisomorphic irreducible modules and
$M_i$ are their multiplicity spaces. Schur's lemma gives the following
criterion.
\begin{proposition}\label{prop:reducible-symplectic-criterion}
Let $V$ be an $\s$-module. For each irreducible
$\s$-module $U$, let $m_U$ denote its multiplicity in $V$, with
$m_U=0$ when it does not occur.  Then $V$ is symplectic if and only if
the following conditions hold.
\begin{enumerate}[label=\textup{(\alph*)}]
\item Every irreducible module $U$ with $\nu(U)=1$ occurs in $V$ with even
multiplicity.
\item For every non-self-dual irreducible $U$, one has $m_U=m_{U^*}$.
\end{enumerate}
There is no parity restriction on the multiplicity of an irreducible module
with indicator $-1$.
\end{proposition}

\begin{proof}
Write an isotypic component corresponding to a self-dual irreducible module
$U$ as $U\otimes M$, where $M$ is the multiplicity space.  By Schur's lemma,
every invariant bilinear form on this component is the tensor product of the
invariant form on $U$ with a bilinear form on $M$.  The tensor product is
alternating precisely when one factor is alternating and the other is
symmetric.  Thus a symplectic $U$ requires a nondegenerate symmetric form on
$M$, which exists in every dimension, whereas an orthogonal $U$ requires a
nondegenerate alternating form on $M$, which exists precisely when
$\dim M$ is even.

If $U$ is not self-dual, an invariant form pairs its isotypic component only
with the component corresponding to $U^*$.  It is nondegenerate precisely
when the two multiplicity spaces have the same dimension.  In that case a
perfect pairing between them gives an invariant alternating form on the sum
of the two components.  Different duality orbits are mutually orthogonal,
and the result follows.
\end{proof}

Thus a symplectic module need not be a sum of irreducible symplectic
modules: orthogonal components can occur with even multiplicity, and a
non-self-dual component can be paired with its dual. For example,
$W\oplus W^*$ has the canonical symplectic form
\begin{equation}\label{eq:canonical-symplectic-module}
\omega((u,\varphi),(v,\psi))=\psi(u)-\varphi(v).
\end{equation}
By contrast, a \emph{multiplicity-free sum of irreducible symplectic
modules} means a sum of pairwise nonisomorphic irreducible modules,
each of which is itself symplectic. This is the hypothesis of our
Heisenberg rigidity theorems.
\begin{lemma}\label{lem:symplectic-not-adjoint}
If $V$ is an irreducible symplectic $\s$-module, then
\[
\Hom_\s(V,\s)=0.
\]
\end{lemma}

\begin{proof}
Suppose that there exists a nonzero $\s$-homomorphism from $V$ to $\s$.
Since $V$ is irreducible, it identifies $V$ with an irreducible submodule of
the adjoint module $\s$.  Such a submodule is a simple ideal of $\s$ and
therefore carries a nondegenerate invariant symmetric bilinear form, given
by the restriction of its Killing form.  On the other hand, $V$ carries a
nondegenerate invariant alternating form.  This is impossible: by Schur's
lemma, the space $\Hom_\s(V,V^*)$ is one-dimensional, so two nonzero
invariant bilinear forms on $V$ must be proportional.
\end{proof}
\begin{lemma}\label{lem:symplectic-parity}
If $V$ is a direct sum of irreducible symplectic $\s$-modules, then
\[
\Hom_\s(\Lambda^2V,V)=0.
\]
\end{lemma}
\begin{proof}
Put $h_0=2\rho^\vee$.  By Tits's formula, the highest weight of
each summand takes an odd integer value on $h_0$.  Every root takes
an even integer value on $h_0$, since $\alpha_i(\rho^\vee)=1$
for each simple root.  All weights of $V$ therefore take odd values
on $h_0$, whereas all weights of $\Lambda^2V$ take even values.
An equivariant map preserves these eigenspaces, so it must be zero.
\end{proof}
We record a few families of irreducible symplectic modules.
\begin{enumerate}[label=\textup{(\alph*)}]
\item For $\s=\mathfrak{sl}_2$, the irreducible module $V(n)$ of highest
weight $n$ satisfies
\[
\nu(V(n))=(-1)^n.
\]
Consequently, $V(n)$ is symplectic precisely when $n$ is odd, or
equivalently, when $\dim V(n)=n+1$ is even.

\item The defining module $\C^{2n}$ of $\mathfrak{sp}_{2n}$ is irreducible
and symplectic for every $n\geq1$.

\item For $\s=\mathfrak{sl}_{4k+2}$, with $k\geq0$, the irreducible module
\[
\Lambda^{2k+1}\C^{4k+2}
\]
is symplectic.  Its invariant form is induced by the exterior product
\[
\Lambda^{2k+1}\C^{4k+2}\times
\Lambda^{2k+1}\C^{4k+2}
\longrightarrow
\Lambda^{4k+2}\C^{4k+2}\simeq\C,
\]
which is alternating because $2k+1$ is odd.
\end{enumerate}

\subsection{An equivariant Jacobi map}

For an $\s$-module $V$, define
\begin{equation}\label{eq:J-map}
\begin{split}
J:\Hom_\s(\Lambda^2V,\s)&\longrightarrow\Hom_\s(\Lambda^3V,V),\\
J(\varphi)(u,v,w)&=\varphi(u,v)\cdot w+
\varphi(v,w)\cdot u+\varphi(w,u)\cdot v.
\end{split}
\end{equation}
The equation $J(\varphi)=0$ is the Jacobi identity on $V$ for the
bracket on $\s\oplus V$ obtained by setting $[u,v]=\varphi(u,v)$.

\begin{lemma}[Killing-form criterion]\label{lem:Jacobi-injective}
If $V$ is a multiplicity-free sum of irreducible symplectic
$\s$-modules, then the map $J$ in \eqref{eq:J-map} is injective.
\end{lemma}
\begin{proof}
Let $\varphi\in\Hom_\s(\Lambda^2V,\s)$ satisfy $J(\varphi)=0$.
On $\mathfrak l=\s\oplus V$, retain the bracket of $\s$ and its action
on $V$, and put $[u,v]=\varphi(u,v)$. Equivariance of $\varphi$ gives
Jacobi on one input in $\s$ and two in $V$, while $J(\varphi)=0$ gives
Jacobi on three inputs in $V$. Hence $\mathfrak l$ is a Lie algebra.

Let $B$ be its Killing form and write $V=\bigoplus_iV_i$ with the
$V_i$ pairwise nonisomorphic irreducible symplectic modules. There is no
nonzero invariant symmetric bilinear form on $V$: its restriction to
each $V_i$ must vanish by Schur's lemma, and its cross terms vanish
because $V_i\not\simeq V_j^*\simeq V_j$ for $i\neq j$. Thus
$B(V,V)=0$. Also $B(V,\s)=0$, by
Lemma~\ref{lem:symplectic-not-adjoint} and the self-duality of the
adjoint module. Therefore $V$ is contained in the ideal
$\mathfrak k=\ker B$.

The ideal $\mathfrak k$ is solvable. Indeed, $\operatorname{ad}(x)$
acts trivially on $\mathfrak l/\mathfrak k$ for $x\in\mathfrak k$,
so the Killing form of $\mathfrak k$ is the restriction of $B$ and is
zero. Cartan's solvability criterion applies. Consequently
$\mathfrak k\cap\s$ is a solvable ideal of $\s$ and must be zero.
It follows that
\[
\operatorname{im}\varphi=[V,V]\subseteq\mathfrak k\cap\s=0.
\]
\end{proof}

\subsection{The common deformation source}

Let $(V,\omega)$ be a nonzero symplectic $\s$-module. Put
\[
\mathfrak r_0=V\oplus\C a,\qquad
\g_0=(\s\ltimes V)\times\C a=\s\ltimes\mathfrak r_0,
\]
where $\mathfrak r_0$ is abelian. Define alternating maps on $\g_0$ by
\[
\tau(a,v)=v,\qquad \eta(u,v)=\omega(u,v)a,
\]
and set them equal to zero on all unspecified pairs.
\begin{lemma}\label{lem:basic-deformations}
Let $\mu_0$ be the Lie bracket of $\g_0$. For every $t\in\C$, the bilinear
maps
\[
\mu_0+t\tau
\qquad\text{and}\qquad
\mu_0+t\eta
\]
are Lie brackets.
\end{lemma}

\begin{proof}
For the first family, the only new bracket is
\[
[a,v]_t=t v.
\]
The endomorphism $t\operatorname{Id}_V$ commutes with the action of $\s$ on
$V$, and hence $a$ acts as a derivation of the abelian Lie algebra $V$ in a
way compatible with the action of $\s$. Thus the Jacobi identity holds.

For the second family, the only new bracket is
\[
[v,v']_t=t\omega(v,v')a.
\]
For this bracket, $[a,\s\oplus V\oplus\C a]_t=0$. Hence the
Jacobi identity for three elements of $V\oplus\C a$ is automatic. The Jacobi identity involving an element $x\in\s$ and two
elements $v,v'\in V$ is precisely the identity
\[
\omega(x\cdot v,v')+\omega(v,x\cdot v')=0.
\]
This follows from the $\s$-invariance of $\omega$.
\end{proof}

For $t\neq0$, an isomorphism from the toral member at parameter $t$
to that at parameter one fixes $\s\oplus V$ and sends $a$ to $ta$.
For the Heisenberg family it sends $a$ to $t^{-1}a$ instead. Thus the
nonzero targets are respectively $\g_\tau$ and $\g_\eta$.

\begin{remark}[The Grunewald--O'Halloran construction]
\label{rem:Grunewald-OHalloran}
If $I$ is a codimension-one ideal of $(\g,\mu)$, write
$\g=\C x\oplus I$ and take $D\in\Der(I)$. The alternating map
\[
\varphi_D(x,y)=D(y),\qquad \varphi_D(I,I)=0
\]
defines a linear deformation $\mu+t\varphi_D$: it replaces the
derivation $\operatorname{ad}_\mu(x)|_I$ by
$\operatorname{ad}_\mu(x)|_I+tD$
\cite{GrunewaldOHalloran1993,HerreraGranadaTirao2016}.
For the source above, take $I=\s\ltimes V$, $x=a$, and
$D|_\s=0$, $D|_V=\operatorname{Id}_V$. This gives $\tau$ and the
original action of $a$ on $I$ is zero. The map $\eta$ instead changes
the bracket on $\Lambda^2V$ and takes values in the complementary line.
\end{remark}

\begin{remark}
The two directions do not give a two-parameter linear family. The
Jacobi expression for $\mu_0+\lambda\tau+t\eta$ on $(a,u,v)$ is
$-2\lambda t\omega(u,v)a$, so this is a Lie bracket precisely when
$\lambda t=0$.
\end{remark}

\section{Toral extensions}\label{sec:toral}

\subsection{A single toral generator}

Let $V$ be any nonzero $\s$-module and put
$A=\operatorname{End}_\s(V)$. No symplectic assumption is needed in
this section.
\begin{proposition}\label{prop:reducible-tau}
Let $V$ be a nonzero $\s$-module and let
\[
\g_\tau=(\s\oplus\C a)\ltimes V,
\qquad [a,v]=v.
\]
Then
\begin{equation}\label{eq:reducible-tau-cohomology}
H^0(\g_\tau,\g_\tau)=0,
\qquad
H^1(\g_\tau,\g_\tau)\simeq A/\C\operatorname{Id}_V,
\qquad
H^2(\g_\tau,\g_\tau)\simeq A/\C\operatorname{Id}_V.
\end{equation}
If
\begin{equation}\label{eq:reducible-isotypic-decomposition}
V=\bigoplus_{i=1}^k U_i^{\oplus m_i},
\end{equation}
where the $U_i$ are pairwise nonisomorphic irreducible modules, then
\begin{equation}\label{eq:reducible-tau-dimensions}
\dim H^1(\g_\tau,\g_\tau)
=\dim H^2(\g_\tau,\g_\tau)
=\sum_{i=1}^k m_i^2-1.
\end{equation}
\end{proposition}

\begin{proof}
Put $\mathfrak r_\tau=\C a\ltimes V$. Corollary~\ref{cor:HS-low-degree}
gives
\[
H^q(\g_\tau,\g_\tau)
\simeq H^q(\mathfrak r_\tau,\g_\tau)^\s,
\qquad q=0,1,2.
\]
Cartan's homotopy formula \eqref{eq:Cartan-homotopy}
reduces the calculation to
\[
C^q_0=\{f\in\Hom_\s(\Lambda^q\mathfrak r_\tau,\g_\tau):
\mathcal L_a f=0\}.
\]  In degrees zero through
three its terms are
\begin{align*}
C^0_0&=\C a,\\
C^1_0&=\C\varepsilon\oplus A,\\
C^2_0&=\{\delta_T:T\in A\}\simeq A,\\
C^3_0&=0.
\end{align*}
Here all unspecified values are zero, and
\[
\varepsilon(a)=a,
\qquad
T(v)=Tv,
\qquad
\delta_T(a,v)=Tv.
\]
Indeed, a zero-weight one-cochain can only send $a$ to $\C a$ and $V$ to
$V$, while a zero-weight two-cochain can only send $a\wedge V$ to $V$.
Equivariance identifies the latter two types of maps with $A$.  A threefold exterior product is either in $\Lambda^3V$ (weight
three) or in $a\wedge\Lambda^2V$ (weight two); a product with two
copies of $a$ is zero.  Since the target has only weights zero and one,
$C^3_0=0$.

For a one-cochain $f$, we use
\[
(df)(x,y)=[x,f(y)]-[y,f(x)]-f([x,y]).
\]
A direct evaluation on $v\in V$ and on $(a,v)$ gives
\begin{equation}\label{eq:reducible-tau-differential}
d a=-\operatorname{Id}_V,
\qquad
d(\lambda\varepsilon+T)
=\lambda\delta_{\operatorname{Id}_V},
\qquad
d\delta_T=0.
\end{equation}
Consequently the map $C^0_0\to C^1_0$ is injective, the kernel in degree one
is $A$, and the space of boundaries is $\C\operatorname{Id}_V$.  In degree two every $\delta_T$ is a cocycle, while
the boundaries are exactly the scalar multiples of
$\delta_{\operatorname{Id}_V}$. This proves
\eqref{eq:reducible-tau-cohomology}.
Finally, $A\simeq\bigoplus_i M_{m_i}(\C)$ for the decomposition
\eqref{eq:reducible-isotypic-decomposition}, which gives
\eqref{eq:reducible-tau-dimensions}.  The subtraction of one comes
from the scalar identity: it is a one-coboundary, and its corresponding
two-cochain is $d\varepsilon$.  This remains true when $V$ contains
trivial one-dimensional summands.
\end{proof}

\begin{corollary}\label{cor:toral-rigidity}
The algebra $\g_\tau$ is algebraically rigid if and only if $V$ is
irreducible. In that case $\g_\tau$ is also complete.
\end{corollary}
\begin{proof}
By complete reducibility, $\dim A=\sum_i m_i^2$ equals one precisely
when $V$ is irreducible. Apply Proposition~\ref{prop:reducible-tau}.
\end{proof}

Every class represented by $\delta_T$ integrates by replacing the action
of $a$ with $\operatorname{Id}_V+tT$.
\begin{corollary}\label{cor:reducible-tau-not-rigid}
If $V$ is reducible, then $\g_\tau$ is not geometrically rigid.
\end{corollary}

\begingroup
\setlength{\emergencystretch}{2em}
\begin{proof}
Choose a proper nonzero submodule and an $\s$-equivariant projection $P$ onto
it.  For $t\in\C$, replace the action of $a$ by
\begin{equation}\label{eq:reducible-tau-deformation}
[a,v]_t=(\operatorname{Id}_V+tP)v.
\end{equation}
Because $P$ commutes with $\s$, this defines a linear deformation whose
infinitesimal part is $\delta_P$. For $t\neq-1$, this action is
invertible.  The solvable radical is
$\C a\ltimes V$, and any nilpotent ideal is contained in it.  If such
an ideal contained $ca+v$ with $c\neq0$, its brackets with $V$ would
contain $V$, on which $\operatorname{ad}(ca+v)$ is invertible.  This
contradicts nilpotence.  Hence the nilradical is $V$ and is characteristic.  The center of the quotient by $V$ is the
line spanned by the class of $a$. Hence an isomorphism with $\g_\tau$,
the member at $t=0$ of this new family, would have to conjugate a nonzero
scalar multiple of
$\operatorname{Id}_V+tP$ to $\operatorname{Id}_V$.  For $t\neq0,-1$, this is impossible: the action has two distinct
eigenvalues, $1$ and $1+t$, on nonzero subspaces.  Thus arbitrarily close members of the
family are not isomorphic to $\g_\tau$, and $\g_\tau$ is not geometrically rigid.
\end{proof}
\endgroup

\subsection{Independent toral generators}

Choose a decomposition $V=\bigoplus_{i=1}^rV_i$ into irreducible
components, allowing repetitions, and let
$Z_r=\bigoplus_{i=1}^r\C a_i$ be a trivial $\s$-module. The
projections onto the $V_i$ define commuting semisimple derivations of
$\s\ltimes V$. Hence
\begin{equation}\label{eq:full-torus-definition}
\g_\tau^{(r)}=(\s\oplus Z_r)\ltimes V,
\qquad [a_i,v_j]=\delta_{ij}v_j
\end{equation}
is a Lie algebra. It is the nonzero target of the linear deformation
of $(\s\ltimes V)\times Z_r$ defined by
$\tau^{(r)}(a_i,v_j)=\delta_{ij}v_j$.

An \emph{external torus} is an abelian subalgebra of semisimple
derivations with zero intersection with the inner derivations. The
derivations just described span a maximal torus in
$A\simeq\bigoplus_i\operatorname{End}(M_i)$ and form an external torus
of $\s\ltimes V$. Indeed, an inner derivation vanishing on $\s$ is
defined by a vector in $V^\s$ and acts trivially on $V$. The scalar
derivation defining $\g_\tau$ is the sum of these independent
generators; its span is maximal precisely when $V$ is irreducible.
This identifies the constructions within the external-torus framework
of \cite{Carles1984,Carles1985}. The following direct calculation gives
the cohomology of the full-torus extension; compare also the broader
cohomological setting of \cite{Carles1993}.
\begin{proposition}[The full-torus extension]\label{prop:full-torus}
Let $V=\bigoplus_{i=1}^rV_i$ be any decomposition into nonzero
irreducible $\s$-modules, allowing isomorphic or trivial summands.
Let $Z_r$ act by $[a_i,v_j]=\delta_{ij}v_j$, and put
$\g_\tau^{(r)}=(\s\oplus Z_r)\ltimes V$. Then
\[
H^q(\g_\tau^{(r)},\g_\tau^{(r)})=0\qquad(q\geq0).
\]
In particular, this extension is complete and algebraically rigid.
\end{proposition}
\begin{proof}
Apply Theorem~\ref{thm:HS-factorization} with semisimple factor $\s$
and radical $\mathfrak r_\tau^{(r)}=Z_r\ltimes V$.
Cartan's homotopy formula for the commuting operators $a_i$ reduces
$\bigl(C^\bullet(\mathfrak r_\tau^{(r)},\g_\tau^{(r)})\bigr)^\s$ to its common zero-weight part.
Let $\alpha_1,\ldots,\alpha_r\in Z_r^*$ be the dual basis,
so $\alpha_i(a_j)=\delta_{ij}$.  The output weights are zero
and the $\alpha_i$; hence a cochain of weight zero has either only
toral inputs and toral output, or exactly one input in $V_i$ and
output in that same $V_i$.  Even if two modules are isomorphic,
their torus weights are different.  Schur's lemma gives
\[
C^q_0\simeq
(\Lambda^q Z_r^*\otimes Z_r)
\oplus(\Lambda^{q-1} Z_r^*\otimes\C^r),
\]
where $\Lambda^{-1} Z_r^*=0$.
If $F\colon\Lambda^q Z_r\to Z_r$ is in the
first summand, then
\[
(dF)(z_1,\ldots,z_q,v_i)
=(-1)^{q+1}\alpha_i(F(z_1,\ldots,z_q))v_i.
\]
The differential is zero on the second summand: the action on its
output cancels the action on its one $V_i$ input, and $V$ is abelian.
Thus the differential identifies the first summand in degree $q$
isomorphically with the second summand in degree $q+1$.
The whole zero-weight complex is acyclic.  The factorization theorem
then gives the assertion in every degree.
\end{proof}

\section{Heisenberg extensions}\label{sec:heisenberg}

\subsection{The general kernel formula}

Let $(V,\omega)$ be a nonzero symplectic $\s$-module and put
$\g_\eta=\s\ltimes\h(V,\omega)$, with
$[u,v]=\omega(u,v)a$. We use the spaces
\begin{equation}\label{eq:reducible-eta-spaces}
A=\operatorname{End}_\s(V),\qquad
\Omega_\s(V)=\Hom_\s(\Lambda^2V,\C),\qquad
S=\Hom_\s(\Lambda^2V,\s).
\end{equation}
The symplectic centralizer of the $\s$-action is
\begin{equation}\label{eq:symplectic-centralizer}
\mathfrak c_\omega
=\{T\in A:\omega(Tu,v)+\omega(u,Tv)=0\text{ for all }u,v\in V\}.
\end{equation}
Define $\mathcal B:A\to\Omega_\s(V)$ by
\begin{equation}\label{eq:L-map}
\mathcal B(T)(u,v)=\omega(Tu,v)+\omega(u,Tv).
\end{equation}
This map is surjective. Indeed, for $\beta\in\Omega_\s(V)$ write
$\beta(u,v)=\omega(Ru,v)$. Invariance gives $R\in A$, and alternation
gives $\omega(Ru,v)=\omega(u,Rv)$, so $\mathcal B(R)=2\beta$.
The dilation derivation is given by
\begin{equation}\label{eq:eta-outer-derivation}
D|_\s=0,\qquad D(v)=v,\qquad D(a)=2a.
\end{equation}
\begin{proposition}\label{prop:reducible-eta}
Suppose that $V^\s=0$ and $\Hom_\s(V,\s)=0$.
For $\g_\eta=\s\ltimes\h(V,\omega)$ one has
\begin{equation}\label{eq:reducible-eta-cohomology}
H^0(\g_\eta,\g_\eta)\simeq\C,\qquad
H^1(\g_\eta,\g_\eta)\simeq\C[D]\oplus\mathfrak c_\omega.
\end{equation}
Here $D$ is given by \eqref{eq:eta-outer-derivation}, and $T\in\mathfrak c_\omega$
represents the derivation
\begin{equation}\label{eq:centralizer-derivations}
D_T|_\s=0,\qquad D_T(v)=Tv,\qquad D_T(a)=0.
\end{equation}
If, in addition, $\Hom_\s(\Lambda^2V,V)=0$, define
\begin{equation}\label{eq:K-map}
K(T)(u,v,w)=\omega(u,v)Tw+\omega(v,w)Tu+\omega(w,u)Tv.
\end{equation}
Let $J$ be the map in \eqref{eq:J-map}.
Then
\begin{equation}\label{eq:eta-kernel-formula}
H^2(\g_\eta,\g_\eta)\simeq
\ker\bigl(K+J:\mathfrak c_\omega\oplus S
                 \longrightarrow\Hom_\s(\Lambda^3V,V)\bigr).
\end{equation}
\end{proposition}

\begin{proof}
Put $\mathfrak r_\eta=\h(V,\omega)=V\oplus\C a$, with
$[u,v]=\omega(u,v)a$, so that $\g_\eta=\s\ltimes\mathfrak r_\eta$.
Corollary~\ref{cor:HS-low-degree} reduces the calculation to the
invariant complex $\bigl(C^\bullet(\mathfrak r_\eta,\g_\eta)\bigr)^\s$.
In degrees zero and one it has terms
\[
\bigl(C^0(\mathfrak r_\eta,\g_\eta)\bigr)^\s=\C a,\qquad
\bigl(C^1(\mathfrak r_\eta,\g_\eta)\bigr)^\s=A\oplus\C q_a,
\]
where $T\in A$ is extended by $T(a)=0$, and $q_a(a)=a$, $q_a(V)=0$.
With $\eta_\beta(u,v)=\beta(u,v)a$ and all other values zero,
\begin{equation}\label{eq:reducible-eta-first-differential}
dT=\eta_{\mathcal B(T)},\qquad dq_a=-\eta_\omega.
\end{equation}
Thus a one-cochain $(T,c)$ is a cocycle precisely when
$\mathcal B(T)=c\omega$.  Since $\mathcal B(\operatorname{Id}_V)=2\omega$,
it decomposes uniquely as
\[
(T,c)=\frac c2(\operatorname{Id}_V,2)+(T-\tfrac c2\operatorname{Id}_V,0),
\]
with the second term in $\mathfrak c_\omega$.
The invariant zero-cochain $a$ has zero differential because
$[a,\g_\eta]=0$.
This proves the assertions in degrees zero and one.

Under the additional hypothesis, every invariant two-cochain is uniquely
\[
\eta_\beta+\delta_T+\varphi,\qquad
\beta\in\Omega_\s(V),\ T\in A,\ \varphi\in S,
\]
where $\delta_T(a,v)=Tv$ and $\varphi$ is zero on $a\wedge V$.
Every $\eta_\beta$ is a cocycle and, by surjectivity of $\mathcal B$,
a coboundary.  The remaining differential is
\begin{align}
d(\delta_T+\varphi)(a,u,v)
&=-\mathcal B(T)(u,v)a,\label{eq:delta-T-a-v-v'}\\
d(\delta_T+\varphi)(u,v,v')
&=-K(T)(u,v,v')-J(\varphi)(u,v,v').
\label{eq:delta-T-u-v-v'}
\end{align}
These are all possible types of input triples.  Therefore the cocycle
conditions are $T\in\mathfrak c_\omega$ and $K(T)+J(\varphi)=0$.
The invariant two-coboundaries have only the $\Omega_\s(V)$ component,
which proves \eqref{eq:eta-kernel-formula}.
\end{proof}
For $\dim V\geq4$, $K$ is injective on $\operatorname{End}(V)$.
For $v'\neq0$, the quotient $v'^\perp/\C v'$ is symplectic of dimension
$\dim V-2$, so there are $u,v\in v'^\perp$ with $\omega(u,v)=1$.
Then $K(T)(u,v,v')=Tv'$.  Nevertheless, a nonzero $K(T)$ may cancel
$J(\varphi)$ in \eqref{eq:eta-kernel-formula}; injectivity of $K$
alone is not a vanishing argument when $S\neq0$.

\subsection{Multiplicity-free rigidity}

\begin{theorem}\label{thm:multiplicity-free-eta}
Let $V=V_1\oplus\cdots\oplus V_r$ be a nonzero multiplicity-free sum
of irreducible symplectic $\s$-modules. For every invariant symplectic
form $\omega$ on $V$, the algebra
$\g_\eta=\s\ltimes\h(V,\omega)$ satisfies
\begin{equation}\label{eq:multiplicity-free-eta-cohomology}
H^0(\g_\eta,\g_\eta)=\C a,\qquad
H^1(\g_\eta,\g_\eta)=\C[D],\qquad
H^2(\g_\eta,\g_\eta)=0.
\end{equation}
In particular, $\g_\eta$ is algebraically and geometrically rigid.
\end{theorem}
\begin{proof}
Each $V_i$ is nontrivial, and Lemma~\ref{lem:symplectic-not-adjoint}
gives $\Hom_\s(V,\s)=0$. Lemma~\ref{lem:symplectic-parity} gives
$\Hom_\s(\Lambda^2V,V)=0$. By Schur's lemma,
\[
A=\bigoplus_{i=1}^r\C\operatorname{Id}_{V_i}.
\]
There are no invariant cross pairings between different $V_i$, so
$\omega$ restricts nondegenerately to every $V_i$. If
$T=\bigoplus_i t_i\operatorname{Id}_{V_i}$ belongs to
$\mathfrak c_\omega$, then $2t_i\omega|_{V_i}=0$ for every $i$.
Thus $\mathfrak c_\omega=0$.

Proposition~\ref{prop:reducible-eta} now gives the first two assertions
and identifies $H^2(\g_\eta,\g_\eta)$ with $\ker J$.
Lemma~\ref{lem:Jacobi-injective} proves the remaining vanishing.
\end{proof}

For an irreducible symplectic module, both nonzero targets of
Lemma~\ref{lem:basic-deformations} are therefore algebraically rigid.
For a multiplicity-free sum of two or more irreducible symplectic
modules, the Heisenberg target is algebraically rigid whereas the
single-generator toral target is geometrically nonrigid.

\begin{remark}\label{rem:perfect-heisenberg-target}
Whenever $V^\s=0$, the algebra $\g_\eta$ is perfect, and its radical
and nilradical are both $\h(V,\omega)$. Indeed,
$[\s,V]=V$ by complete reducibility, and $[V,V]=\C a$ by
nondegeneracy. Every nilpotent ideal maps to zero in the semisimple
quotient, so it lies in the Heisenberg radical. Thus these targets
belong to the perfect case in Carles's structural alternatives for
geometrically rigid algebras \cite{Carles1984}.
\end{remark}

\begin{remark}[The stronger adjoint exclusion]\label{rem:eta-exclusion-scope}
The condition $V^\s=0$ and $S=0$ implies the hypotheses of
Theorem~\ref{thm:multiplicity-free-eta}, but is not needed there.
To see this implication, let $\kappa$ be the Killing form of $\s$.
For $T\in\mathfrak c_\omega$, the form
$b_T(u,v)=\omega(Tu,v)$ is symmetric and invariant. Define
$\Phi_T:\Lambda^2V\to\s$ by
\[
\kappa(\Phi_T(u,v),x)=b_T(x\cdot u,v).
\]
Symmetry and invariance make $\Phi_T$ alternating and equivariant.
If $\Phi_T=0$, then $T(x\cdot u)=0$ for all $x,u$; since
$[\s,V]=V$, one has $T=0$. Hence $\mathfrak c_\omega$ injects into
$S$, and $S=0$ forces $\mathfrak c_\omega=0$.

The multiplicity-space description in
Proposition~\ref{prop:reducible-symplectic-criterion} identifies the
factors of $\mathfrak c_\omega$ with $\mathfrak{so}(M,b)$ for a
symplectic irreducible component, $\mathfrak{sp}(M,b)$ for an
orthogonal one, and $\mathfrak{gl}(M)$ for a non-self-dual pair. Here
$b$ is the corresponding nondegenerate symmetric or alternating form
on the multiplicity space. These factors all vanish precisely when
every component is symplectic and has multiplicity one. In particular,
the stronger exclusion $S=0$ implies this condition. Examples with
$S\neq0$ are given in Section~\ref{sec:examples}.
\end{remark}

\subsection{Separate Heisenberg centers}

Suppose $V=\bigoplus_{i=1}^rV_i$, with each $V_i$ a nonzero
symplectic $\s$-module carrying a form $\omega_i$. Let
$Z_r=\bigoplus_{i=1}^r\C a_i$ be a trivial module and put
\begin{equation}\label{eq:blockwise-heisenberg-target}
\g_\eta^{(r)}=\s\ltimes\bigoplus_{i=1}^r\h(V_i,\omega_i),
\qquad [u_i,v_i]=\omega_i(u_i,v_i)a_i.
\end{equation}
Different radical blocks commute. The same Levi factor acts on all
blocks; this is not a direct product of $r$ semidirect products.
Its radical is $\mathfrak r_\eta^{(r)}=V\oplus Z_r$.

The algebra \eqref{eq:blockwise-heisenberg-target} is a nonzero target
of the linear deformation of $(\s\ltimes V)\times Z_r$ given by
\[
\eta^{(r)}(u_i,v_i)=\omega_i(u_i,v_i)a_i,
\qquad \eta^{(r)}(V_i,V_j)=0\quad(i\neq j).
\]
Jacobi follows from the centrality of $Z_r$ and the invariance of each
$\omega_i$. Rescaling the $a_i$ identifies all nonzero members.
When the $V_i$ are irreducible, this source also gives the full-torus
target of \eqref{eq:full-torus-definition}.

\begin{theorem}\label{thm:multi-heisenberg}
If $V_1,\ldots,V_r$ are pairwise nonisomorphic irreducible symplectic
$\s$-modules, then
\begin{equation}\label{eq:multi-heisenberg-cohomology}
\begin{split}
H^0(\g_\eta^{(r)},\g_\eta^{(r)})&=Z_r,\\
H^1(\g_\eta^{(r)},\g_\eta^{(r)})&\simeq\C^r,\\
H^2(\g_\eta^{(r)},\g_\eta^{(r)})&=0.
\end{split}
\end{equation}
The first cohomology is generated by the classes of the derivations
$D_i$ defined by $D_i|_{V_i}=\operatorname{Id}_{V_i}$,
$D_i(a_i)=2a_i$, and zero on the remaining components.
In particular, $\g_\eta^{(r)}$ is algebraically rigid.
\end{theorem}
\begin{proof}
Corollary~\ref{cor:HS-low-degree} reduces the calculation to the
$\s$-invariant cochain complex of $\mathfrak r_\eta^{(r)}$.
Schur's lemma gives
\[
\operatorname{End}_\s(V)=\bigoplus_i\C\operatorname{Id}_{V_i},
\qquad \Omega_\s(V)=\bigoplus_i\C\omega_i.
\]
Also $V^\s=0$, $\Hom_\s(V,\s)=0$, and
$\Hom_\s(\Lambda^2V,V)=0$, by the preceding lemmas.

The invariant zero-cochains are $Z_r$ and have zero differential.
An invariant one-cochain has the form
$f(v_i)=p_iv_i$, $f(z)=Qz$, with $p_i\in\C$ and
$Q\in\operatorname{End}(Z_r)$. Its differential is
\begin{equation}\label{eq:multi-heisenberg-first-differential}
(df)(u_i,v_i)=\omega_i(u_i,v_i)(2p_i a_i-Q(a_i)),
\end{equation}
and is zero on the other pairs. Thus it is a cocycle precisely when
$Q(a_i)=2p_i a_i$ for all $i$. This proves the assertions in degrees
zero and one.

Every invariant two-cochain is uniquely a sum
\[
\beta+\delta_\lambda+\gamma+\varphi,
\]
where $\beta(u_i,v_i)=\omega_i(u_i,v_i)b_i$ for arbitrary $b_i\in Z_r$,
$\delta_\lambda(z,v_i)=\lambda_i(z)v_i$ for $\lambda_i\in Z_r^*$,
$\gamma\in\Hom(\Lambda^2Z_r,Z_r)$, and
$\varphi\in\Hom_\s(\Lambda^2V,\s)$. Each map is extended by zero
on the other types of pairs. The stated vanishing results exclude
all further cochain types.

The differential on $(z,u_i,v_i)$ is
\[
-\omega_i(u_i,v_i)\bigl(2\lambda_i(z)a_i+\gamma(a_i,z)\bigr).
\]
Thus a cocycle satisfies $\gamma(a_i,z)=-2\lambda_i(z)a_i$.
Taking $z=a_i$ gives $\lambda_i(a_i)=0$. Taking $z=a_j$ and
exchanging $i,j$ gives
\[
-2\lambda_i(a_j)a_i=\gamma(a_i,a_j)=2\lambda_j(a_i)a_j.
\]
For $i\neq j$, linear independence of $a_i,a_j$ forces both
coefficients to vanish. Hence $\lambda_i=0$ for all $i$ and
$\gamma=0$.

Every $\beta$ is a coboundary: in
\eqref{eq:multi-heisenberg-first-differential}, take $p_i=0$ and
$Q(a_i)=-b_i$. The remaining cocycle equation on three inputs in
$V$ is $J(\varphi)=0$. Lemma~\ref{lem:Jacobi-injective} gives
$\varphi=0$, proving the vanishing in degree two.
\end{proof}

\section{Applications and nonrigid examples}\label{sec:examples}

\subsection{Irreducible families}

\begin{corollary}\label{cor:symplectic-oscillator}
Let $V=\C^{2n}$ be the defining module of $\mathfrak{sp}_{2n}$,
where $n\geq1$. Both nonzero targets associated with
$(\mathfrak{sp}_{2n}\ltimes V)\times\C a$ are algebraically rigid.
The Heisenberg target is the symplectic oscillator algebra
$\mathfrak{sp}_{2n}\ltimes\h_{2n+1}$.
\end{corollary}
\begin{proof}
The defining module is irreducible and symplectic. Apply
Corollary~\ref{cor:toral-rigidity} and
Theorem~\ref{thm:multiplicity-free-eta}. For the terminology, see
\cite{LiuZhao2019}.
\end{proof}

\begin{corollary}\label{cor:conformal-galilei}
Let $m$ be a positive odd integer and let $d=1$ or $d\geq3$. Put
\[
\s_d=\begin{cases}
\mathfrak{sl}_2,&d=1,\\
\mathfrak{sl}_2\oplus\mathfrak{so}_d,&d\geq3,
\end{cases}
\qquad
V_d=\begin{cases}
V(m),&d=1,\\
V(m)\otimes\C^d,&d\geq3,
\end{cases}
\]
where $V(m)$ is the irreducible $\mathfrak{sl}_2$-module of highest
weight $m$. Then both associated nonzero targets are algebraically
rigid. The Heisenberg target is the complexified mass central
extension of the conformal Galilei algebra of spin $m/2$; for $m=1$
this is the Schr\"odinger family.
\end{corollary}
\begin{proof}
For $d\geq3$, the module $V_d$ is irreducible, and the tensor product
of the alternating form on $V(m)$ with the symmetric form on $\C^d$
is symplectic. For $d=1$, $V(m)$ itself is irreducible symplectic.
Apply Corollary~\ref{cor:toral-rigidity} and
Theorem~\ref{thm:multiplicity-free-eta}.
\end{proof}

The Heisenberg rigidity assertions in this corollary recover
\cite[Theorems~5.1 and~5.2]{KhudoyberdiyevJumaniyozov2025}, and for
$m=1$ the corresponding results of \cite{Campoamor2017}. In particular,
the case of a single odd irreducible $\mathfrak{sl}_2$-module is
already included in that literature. The factor $\mathfrak{so}_2$ is
abelian, so the spatial dimension $d=2$ is outside the semisimple-Levi
formulation used here.

\subsection{Multiplicity-free modules and adjoint cochains}

For $\mathfrak{sl}_2$, let $V(n)$ denote the irreducible module of
highest weight $n$. The standard decompositions are
\begin{equation}\label{eq:sl2-Clebsch-Gordan}
V(m)\otimes V(n)=\bigoplus_{j=0}^{\min\{m,n\}}V(m+n-2j)
\end{equation}
and
\begin{equation}\label{eq:sl2-exterior-square}
\Lambda^2V(n)=\bigoplus_{j=0}^{\lfloor(n-1)/2\rfloor}V(2n-2-4j).
\end{equation}
The second follows from the first because the flip acts on the
component $V(2n-2j)$ by $(-1)^j$. For $n=0$ the sum is empty.
Also $\nu(V(n))=(-1)^n$, so
$\bigoplus_n V(n)^{\oplus m_n}$ is symplectic precisely when every
even $n$ has even multiplicity $m_n$.

\begin{corollary}\label{cor:sl2-multiplicity-free}
Let $V=\bigoplus_{i=1}^rV(n_i)$, where the $n_i$ are pairwise
distinct positive odd integers. Then
\[
H^0(\g_\tau,\g_\tau)=0,\qquad
H^1(\g_\tau,\g_\tau)\simeq H^2(\g_\tau,\g_\tau)
\simeq\C^{r-1},
\]
whereas $\g_\eta$ satisfies \eqref{eq:multiplicity-free-eta-cohomology}
for every $r\geq1$. The extension with a separate Heisenberg center
for each $V(n_i)$ satisfies \eqref{eq:multi-heisenberg-cohomology}.
\end{corollary}
\begin{proof}
Apply Proposition~\ref{prop:reducible-tau},
Theorem~\ref{thm:multiplicity-free-eta}, and
Theorem~\ref{thm:multi-heisenberg}.
\end{proof}

The adjoint module $V(2)$ occurs in
$V(n_i)\otimes V(n_j)$ for distinct odd $n_i,n_j$ precisely when
$|n_i-n_j|=2$. Thus $\Hom_\s(\Lambda^2V,\s)$ can be nonzero in this
corollary. Such cochains are eliminated by the Jacobi equation in
Lemma~\ref{lem:Jacobi-injective}, rather than by a representation-theoretic
exclusion. The same phenomenon occurs for other semisimple Levi factors.

\begin{example}\label{ex:semisimple-adjoint-cochains}
Take $\s=\mathfrak{sl}_2\oplus\mathfrak{sl}_2$ and
\[
V=(V(1)\otimes V(2))\oplus(V(3)\otimes V(2)).
\]
The two components are irreducible symplectic and nonisomorphic.
Their cross tensor product contains $V(2)\otimes V(0)$, an adjoint
component of $\s$, by \eqref{eq:sl2-Clebsch-Gordan}. Hence
$\Hom_\s(\Lambda^2V,\s)\neq0$, but both the single-center and
separate-center Heisenberg targets are algebraically rigid.
\end{example}

\subsection{Repeated blocks with separate centers}

Repeated irreducible components need not preserve rigidity. We first
consider the construction with one center for each copy.

\begin{proposition}\label{prop:repeated-blocks-nonrigid}
Let $U$ be a nonzero irreducible symplectic $\s$-module and let
\[
\g_\eta^{(m)}=\s\ltimes\bigoplus_{i=1}^m\h(U_i,\omega_U),
\qquad U_i\simeq U.
\]
For $m\geq3$, this algebra is geometrically nonrigid. Moreover,
\[
\dim H^2(\g_\eta^{(m)},\g_\eta^{(m)})
\geq\frac{m(m-1)(m-2)}2.
\]
\end{proposition}
\begin{proof}
Write $u_i$ for the copy of $u\in U$ in $U_i$, and let
$Z_m=\bigoplus_i\C a_i$. Choose three distinct indices, say $1,2,3$,
and define an alternating map by
\begin{equation}\label{eq:repeated-block-deformation}
\beta(u_1,v_2)=\omega_U(u,v)a_3,
\end{equation}
with all other types of pairs zero. If $\mu$ is the original bracket,
then $\mu+t\beta$ is a Lie bracket for every $t$: all radical brackets
still take values in $Z_m$, and the new bracket is $\s$-equivariant.

The cocycle $\beta$ is not a coboundary. By
Corollary~\ref{cor:HS-low-degree} and exactness of taking invariants,
it is enough to check invariant one-cochains on the radical. Every
such cochain has components in $\operatorname{End}_\s(U^{\oplus m})$
and $\operatorname{End}(Z_m)$, since $U$ is nontrivial and
$\Hom_\s(U,\s)=0$. On $U_1\wedge U_2$, its differential takes values
only in $\C a_1\oplus\C a_2$: the original bracket of these two
blocks is zero, and each remaining term brackets against either
$U_1$ or $U_2$. It cannot equal the nonzero $a_3$ component in
\eqref{eq:repeated-block-deformation}.

If the orbit of $\mu$ were open, the curve $\mu+t\beta$ would lie
in that orbit for $t$ in a neighborhood of zero. Its tangent $\beta$
would then belong to the orbit tangent space, which is the space of
two-coboundaries. This proves geometric nonrigidity.

For the dimension estimate, put $M=\C^m$. The central-valued invariant
two-cochains on $U\otimes M$, extended by zero on inputs in $Z_m$,
are all cocycles and have dimension
\[
\dim(S^2M^*\otimes Z_m)=\frac{m^2(m+1)}2.
\]
The invariant one-cochains have dimension $2m^2$. Their kernel
contains the $m$ independent block dilations, so the boundary space
has dimension at most $2m^2-m$. Subtracting gives the stated bound.
\end{proof}

The proposition does not assert nonrigidity for two repeated blocks.
The number of centers also matters, as the next example illustrates.

\subsection{Two repeated modules with a shared center}
\begin{remark}\label{rem:sl2-higher-multiplicities}
For $V=\bigoplus_n V(n)\otimes M_n$, the Lie algebra of equivariant
symplectic endomorphisms is
\[
\mathfrak c_\omega\simeq
\bigoplus_{n\text{ odd}}\mathfrak{so}(M_n)
\oplus\bigoplus_{n\text{ even}}\mathfrak{sp}(M_n),
\]
with the symmetric or alternating multiplicity-space forms from
Proposition~\ref{prop:reducible-symplectic-criterion}.
For a sum of odd components, $V^{\mathfrak{sl}_2}=0$ and
$\Hom_{\mathfrak{sl}_2}(V,\mathfrak{sl}_2)=
\Hom_{\mathfrak{sl}_2}(\Lambda^2V,V)=0$, so
Proposition~\ref{prop:reducible-eta} computes $H^1(\g_\eta,\g_\eta)$
and gives the kernel formula for $H^2(\g_\eta,\g_\eta)$ with arbitrary
multiplicities.
Repeated odd components introduce adjoint-valued cochains through
\[
\Lambda^2(V(n)\otimes M_n)
\simeq(\Lambda^2V(n)\otimes S^2M_n)
\oplus(S^2V(n)\otimes\Lambda^2M_n),
\]
since $V(2)$ occurs in $S^2V(n)$ for every positive odd $n$.
The kernel formula retains both the centralizer and these adjoint-valued
cochains.  Even components require separate checks of all the
hypotheses of Proposition~\ref{prop:reducible-eta}.
\end{remark}
\begin{example}[Two copies of the defining module]\label{ex:two-standard-copies}
Let $U=\C^2$ with its standard symplectic form $\epsilon$, and
let $V_1,V_2$ be two copies of $U$. Write $V=V_1\oplus V_2$,
with form $\omega=\epsilon\oplus\epsilon$, and denote by $u_i\in V_i$
the copy of $u\in U$.
Then
\[
\Lambda^2V\simeq V(2)\oplus V(0)^{\oplus3}.
\]
Use a single shared generator $a$ for both copies, so that
$[u_i,v_i]=\epsilon(u,v)a$ and $[V_1,V_2]=0$. For
$\g_\eta=\mathfrak{sl}_2\ltimes\h(V,\omega)$,
\[
\dim H^0(\g_\eta,\g_\eta)=1,\qquad
\dim H^1(\g_\eta,\g_\eta)=2,\qquad
\dim H^2(\g_\eta,\g_\eta)=1.
\]
Moreover, $\g_\eta$ admits a linear deformation whose nonzero members
are isomorphic to $\mathfrak{sl}_3(\C)$. In particular, $\g_\eta$ is
not geometrically rigid, whereas every nonzero member is algebraically
and geometrically rigid.

To see this explicitly, identify $\mathfrak{sl}_2=\mathfrak{sp}(U)$
and set
\[
\sigma(u,v)x=\epsilon(u,x)v+\epsilon(v,x)u.
\]
Let $T(u_1)=u_2$ and $T(u_2)=-u_1$, and define an alternating map
$\psi$ by
\[
\psi(a,v)=Tv,\qquad
\psi(u_1,v_2)=-\tfrac13\sigma(u,v),\qquad
\psi(V_i,V_i)=\psi(\mathfrak{sl}_2,\g_\eta)=0.
\]
Here $\mathfrak c_\omega=\C T$ and
$\Hom_{\mathfrak{sl}_2}(\Lambda^2V,\mathfrak{sl}_2)$ is spanned by
the cross-pair cochain $\varphi(u_1,v_2)=\sigma(u,v)$.
The two-dimensional symplectic identity
\[
\epsilon(u,v)v'+\epsilon(v,v')u+\epsilon(v',u)v=0
\]
gives, on a triple with two inputs in the first copy,
\[
J(\varphi)(u_1,v_1,v'_2)
=-3\epsilon(u,v)v'_1=3K(T)(u_1,v_1,v'_2).
\]
The same calculation with the copies exchanged gives
$J(\varphi)=3K(T)$ on all triples.
Since $K(T)\neq0$, the kernel in \eqref{eq:eta-kernel-formula}
is precisely the line represented by $\psi=\delta_T-\tfrac13\varphi$.
This proves the cohomology dimensions.

The map $\psi$ is itself a Lie bracket. Indeed,
$\psi(\mathfrak{sl}_2,\g_\eta)=0$ and
$\psi(\Lambda^2V)\subset\mathfrak{sl}_2$, so the Jacobi identity
on three vectors of $V$ is automatic. The only remaining condition is
$\psi(Tu,v)+\psi(u,Tv)=0$, which follows from the displayed formulas.
Here the vector subspace $\mathfrak{sl}_2$ lies in the center for
the bracket $\psi$ alone, not for the original bracket $\mu$.
If $\mu$ denotes the bracket of $\g_\eta$, then $d\psi=0$ shows that
$\mu+t\psi$ is a linear deformation of $\g_\eta$.
For $t\neq0$, put $u_\pm=u_1\mp i u_2$ and view $\epsilon(v,-)$
as a row covector. An isomorphism onto $\mathfrak{sl}_3(\C)$ is given,
for $x\in\mathfrak{sl}_2$, by
\[
x+\lambda a+u_++v_-\longmapsto
\begin{pmatrix}
x+\tfrac{it\lambda}{3}\operatorname{Id}_U & u\\[3pt]
-\tfrac{4it}{3}\epsilon(v,-) & -\tfrac{2it\lambda}{3}
\end{pmatrix}.
\]
Indeed,
$[a,u_\pm]_t=\pm it\,u_\pm$, the brackets within either copy
$\{u_\pm:u\in U\}$ vanish, and
\[
[u_+,v_-]_t=2\epsilon(u,v)a-\frac{2it}{3}\sigma(u,v).
\]
The identity
$u\,\epsilon(v,-)=\tfrac12\bigl(\sigma(u,v)-\epsilon(u,v)\operatorname{Id}_U\bigr)$
identifies this mixed bracket with the matrix commutator; the other
brackets follow from the $\mathfrak{sl}_2$-invariance of $\epsilon$.
Whitehead's second lemma shows that every nonzero member is
algebraically, and hence geometrically, rigid. Its center is zero,
whereas $Z(\g_\eta)=\C a$, so these members are nonisomorphic to
$\g_\eta$. This proves the geometric nonrigidity of $\g_\eta$.
The cohomology dimensions also agree with the corresponding
eight-dimensional entry in \cite[Proposition~5.3]{BurdeWagemann2024}.
\end{example}

\end{document}